\documentclass[10pt]{amsart}
\usepackage[cp1251]{inputenc}
\usepackage[english,russian]{babel}
\usepackage{amsmath}
\usepackage{amssymb}
\usepackage{amsfonts}

\begin{document}
\sloppy
\begin{center}
{\large\bf $(n+1)$-ARY DERIVATIONS OF SEMISIMPLE FILIPPOV ALGEBRAS}\\

\hspace*{8mm}

{\large\bf Ivan Kaygorodov}\\
e-mail: kib@math.nsc.ru

{\it 
Sobolev Inst. of Mathematics\\ 
Novosibirsk, Russia\\}
\end{center}

\

\medskip

\

\begin{center} {\bf Abstract: }\end{center}                                                                    
We defined $(n+1)$-ary derivations of $n$-ary algebras. 
We described $(n+1)$-derivations of 
simple and semisimple finite-dimensional Filippov algebras over algebraically closed field zero characteristic.
We constructed new examples of non-trivial $(n+1)$-ary derivations of Filippov algebras.

\medskip


\medskip

\section{Введение}
Одним из способов обобщения дифференцирований является $\delta$-диф\-фе\-рен\-ци\-рование. 
Под $\delta$-диф\-ференцированием алгебры $A$, при $\delta$ --- фиксированном элементе основного поля,  мы понимаем линейное отображение 
$\phi: A\rightarrow A,$ такое что для произвольных $x,y \in A$ верно
$$\phi(xy)=\delta(\phi(x)y+x\phi(y)).$$
В свое время, $\delta$-дифференцирования изучались в работах 
\cite{Fil}-\cite{kay_nary}, где были описаны $\delta$-дифференцирования 
первичных лиевых \cite{Fil,Fill}, первичных альтернативных и мальцевских \cite{Filll} алгебр, 
простых \cite{kay_lie, kay_lie2} и первичных \cite{Zus} лиевых супералгебр, 
полупростых конечномерных йордановых алгебр \cite{kay,kay_lie2} и супералгебр \cite{kay, kay_lie2,kay_zh,kay_de5}, 
алгебр Филиппова малых размерностей и простых конечномерных алгебр Филиппова \cite{kay_nary}, 
а также простой тернарной алгебры Мальцева $M_8$ \cite{kay_nary}. 
В частности, были построены примеры нетривиальных $\delta$-дифференцирований 
для некоторых алгебр Ли \cite{Fill,Zus, kay_gendelta}, простых йордановых супералгебр \cite{kay_zh, kay_de5} и 
некоторых $n$-арных алгебр Филиппова \cite{kay_nary}.

В тоже время, $\delta$-дифференцирование является частным случаем квазидифференцирования и обобщенного дифференцирования. 
Под обобщенным дифференцированием $D$ мы понимаем такое линейное отображение, что существуют линейные отображения $E$ и $F,$
связанные с $D$ условием, таким что для произвольных $x,y \in A$ верно
$$D(xy)=E(x)y+xF(y).$$
Если вдобавок к этому $E=F,$ то $D$ --- квазидифференцирование. 
Тройки $(D,E,F),$ где $D$ --- обобщенное дифференцирование, а $E,F$ --- связанные с ним линейные отображения, называются тернарными дифференцированиями. 
Квазидифференцирования, обобщенные дифференцирования и тернарные дифференцирования рассматривались в работах \cite{LL}-\cite{kay_nplus}.
В частности, изучались обобщенные и тернарные дифференцирования 
алгебр Ли \cite{LL},
супералгебр Ли \cite{ZhangRY}, 
ассоциативных алгебр \cite{komatsu_nak_03, komatsu_nakajima04}, 
обобщенных алгебр Кэли-Диксона \cite{JG3},
йордановых алгебр \cite{shest}.

Понятие тернарного дифференцирования для бинарной алгебры допускает обобщение на случай $n$-арных алгебр. 
В данном случае, под $(n+1)$-арным дифференцированием $n$-арной алгебры $A$ мы подразумеваем такой набор $(f_0,f_1, \ldots, f_n) \in End(A)^{n+1},$
что для произвольных $x_1, \ldots, x_n \in A$ верно 
$$f_0[x_1, \ldots, x_n]=\sum\limits_{i=1}^{n}[x_1, \ldots, f_i(x_i), \ldots , x_n].$$
Соответственно, для $(n+1)$-арного квазидифференцирования необходимо дополнительно требовать $f_1=f_2=\ldots=f_n.$ 
Ясно, что если $\psi_1, \ldots, \psi_n$ --- набор элементов центроида $n$-арной алгебры $A$, то 
$(\sum \psi_i, \psi_1, \ldots, \psi_n)$ --- является $(n+1)$-арным дифференцированием алгебры $A$.
В тоже время, если $D$ --- дифференцирование $n$-арной алгебры $A$, то набор $(D, D, \ldots, D)$ --- $(n+1)$-арное дифференцирование
алгебры $A$. 
Приведенные два вида $(n+1)$-арных дифференцирований, как и их линейные комбинации, 
мы будем считать тривиальными. 
Наибольший интерес представляет вопрос нахождения $(n+1)$-арных дифференцирований, отличных от тривиальных.
Легко заметить следующее 

\medskip

{\bf  Утверждение.} Пусть $A$ --- (анти)коммутативная $n$-арная алгебра и $(f_0, f_1, \ldots, f_n)$ --- $(n+1)$-арное дифференцирование, 
тогда для любой подставновки $\sigma \in S_n$ верно, что $(f_0, f_{\sigma(1)}, \ldots, f_{\sigma(n)})$ --- $(n+1)$-арное дифференцирование. 

\medskip

{\bf Доказательство.} Доказательство данного факта вытекает из известного утверждения о разложении произвольной подстановки в произведение 
транспозиций и очевидного утверждения леммы если $\sigma$ --- транспозиция. Утверждение доказано.

\medskip

Обозначим пространство дифференцирований, $\delta$-дифференцирований, квазидифференцирований и обобщенных дифференцирований, 
соответственно, через $Der(A), Der_{\delta}(A), QDer(A)$ и $GDer(A).$
Очевидно, что мы имеем цепочку включений 
\begin{eqnarray}\label{vklyach}
Der(A) \subseteq Der_{\delta}(A) \subseteq QDer(A) \subseteq GDer(A) \subseteq  End(A).
\end{eqnarray}
Стоит отметить, что если $n$-арная алгебра $A$ с ненулевым умножением и характеристика поля отлична от $n-1$, 
то первое включение всегда строгое. 
В противном случае, в силу того, что тождественное отображение явлется элементом центроида и $\frac{1}{n}$-дифференцированием,
мы получили бы противоречение с определением умножения в алгебре $A$. Понятно, что в случае бинарных алгебр ограничение на 
характеристику поля является не существенным.
Ясно, что каждая из алгебр $Der(A), Der_{\delta}(A), QDer(A)$ и $GDer(A)$ относительно коммутаторного умножения становится алгеброй Ли. 
Пусть $Ann(GDer(A))$ --- аннулятор алгебры Ли обобщенных дифференцирований алгебры $A$.
Отметим, что $Ann(GDer(A))$ не тривиален. 
Действительно, там лежат отображения вида $\alpha \cdot id,$ где $\alpha$ --- элемент основного поля.
Обозначим $$\Delta(A)=GDer(A)/Ann(GDer(A)).$$ 
Также нас будет интересовать структура алгебры Ли $\Delta(A).$
 $(n+1)$-Арные дифференцирования исследовались в работе \cite{kg_nary}, где были описаны $4$-арные и обобщенные дифференцирования 
тернарной алгебры Мальцева $M_8$. В частности, было показано, что $\Delta(M_8)=B_3.$

Основной целью данной работы, результаты которой аннонсированы в \cite{kay_nplus},
 является исследование обобщенных и $(n+1)$-арных дифференцирований простых конечномерных 
 $n$-арных алгебр Филиппова над алгебраически замкнутым полем характеристики нуль. 
В работе дано полное описание квазидифференцирований, обобщенных дифференцирований и  $(n+1)$-арных дифференцирований данных алгебр.
В итоге, с использованием результатов \cite{kay_nary}, мы имеем, что для 
простой конечномерной алгебры Филиппова $A$ над алгебраически замкнутым полем характеристики нуль цепочка включений (\ref{vklyach}) имеет следующий вид
$$Der(A) \subset Der_{\delta}(A) \subset QDer(A) = GDer(A) = End(A),$$
то есть, в последних двух случаях у нас наблюдается равенство, а остальные все включения строгие.
Как следствие, мы получаем описание $(n+1)$-арных дифференцирований полупростых конечномерных алгебр Филиппова (не являющихся простыми) 
над алгебраически замкнутым полем
характеристики нуль, для которых цепочка включений (\ref{vklyach}) имеет следующий вид
$$Der(A) \subset Der_{\delta}(A) \subset QDer(A) = GDer(A) \subset  End(A).$$

\section{Основные леммы.}

Алгеброй Филиппова, определение которой появляется в \cite{fil_nar}, 
называется алгебра $L$ с одной антикоммутативной $n$-арной операцией $[x_1, \ldots, x_n]$, удовлетворяющей тождеству

$$[[x_1,\ldots, x_n],y_2, \ldots, y_n]=\sum\limits_{i=1}^n[x_1,\ldots, [x_i,y_2,\ldots, y_n], \ldots, x_n].$$

Примером $(n+1)$-мерной $n$-арной алгебры Филиппова является алгебра, которую мы будем обозначать $A_{n+1}.$
Известно, что в $A_{n+1}$ можно выбрать базис 
$$\{e\}=\{e_1, \ldots, e_{n+1} \}$$ 
со следующей таблицей умножения:
$$[e_1, \ldots, \hat{e}_i, \ldots, e_{n+1}]=(-1)^{n+i+1}e_i,$$
где через $\hat{e}_i$ обозначается отсутствие элемента $e_i$ в $n$-арном произведении. 
Как было отмечено в работе \cite{Ling}, алгебрами типа $A_{n+1}$ исчерпываются все простые конечномерные $n$-арные алгебры 
Филиппова над алгебраически замкнутым полем характеристики нуль.

\medskip

{\bf  Лемма 1.}
Пусть $A=(a_{ij}),$ где $a_{ii}=0,$  --- матрица линейного преобразования $f_0$ и 
$B=-A^{T}$ --- матрица линейного преобразования $f$, тогда $(f_0, f, \ldots, f)$ --- $(n+1)$-арное дифференцирование
$n$-арной алгебры $A_{n+1}.$

\medskip

{\bf Доказательство.} Достаточно заметить, что выполнена следующая цепочка равенств
$$(-1)^{n+i+1}f_0[e_1,  \ldots, \hat{e}_i, \ldots, e_{n+1}]
=f_0(e_i)=\sum\limits_{j=1}^{n+1} a_{ji}e_j= $$
$$\sum\limits_{j=1}^{n+1} (-1)^{n+j+1} a_{ji} [e_1, \ldots, \hat{e}_j, \ldots, e_{n+1}]=$$
$$\sum\limits_{j=1}^{n+1} (-1)^{n+i+1} [e_1, \ldots, \underbrace{-\sum\limits_{k=1}^{n+1} a_{jk}e_k}\limits_{j}, \ldots, \hat{e}_i, \ldots, e_{n+1}]=$$
$$(-1)^{n+i+1}\sum\limits_{j=1}^{n+1}[e_1, \ldots, \underbrace{f(e_j)}\limits_{j}, \ldots, \hat{e}_i, \ldots, e_{n+1}].$$

А также, что цепочка равенств 
$$[f(e_i), e_i, e_{i_1}, \ldots, e_{i_{n-2}} ] - [f(e_i), e_i, e_{i_1}, \ldots, e_{i_{n-2}} ]=$$
$$[f(e_i), e_i, e_{i_1}, \ldots, e_{i_{n-2}} ] + [e_i, f(e_i), e_{i_1}, \ldots, e_{i_{n-2}} ]+ \ldots+[e_i, e_i, e_{i_1}, \ldots, f(e_{i_{n-2}}) ]=$$
$$f_0 [e_i, e_i, e_{i_1}, \ldots, e_{i_{n-2}} ]$$
выполняется в силу того, что левое и правое выражения в цепочке равны нулю. 
Ясно, что если среди $\{ i_k\}$ найдется индекс, совпадающий с $i$, то цепочка равенств не нарушается.
Лемма доказана.

\medskip

 Далее, для произвольного вектора $v$ через $v|_{e_i}$ мы будем обозначать коэффициент при векторе $e_i$ в его разложении по базису $\{e\}.$
  
\medskip

{\bf  Лемма 2.} Пусть $(f_0, f_1, \ldots, f_n)$ --- $(n+1)$-арное дифференцирование $n$-арной алгебры $A_{n+1}$, 
тогда 
$$(f_0, f_1, \ldots, f_n)= (g_0, g, \ldots, g) +(f_0^*, f_1^*, \ldots, f_n^*),$$
где $f_j^*(e_i)|_{e_i}=f_j^*(e_i)$ и $g_0(e_i)|_{e_i}=g(e_i)|_{e_i}=0.$

\medskip

{\bf Доказательство.}
Отметим, что если $i_k=i_m=i,$ то
$$0=f_0[e_{i_1}, \ldots, e_{i_k}, \ldots, e_{i_m},\ldots, e_{i_n}]=
\sum[e_{i_1}, \ldots, f_j(e_{i_j}), \ldots, e_{i_n}]=$$
$$[e_{i_1}, \ldots, f_k(e_{i_k}), \ldots, e_{i_m}, \ldots, e_{i_n}] + [e_{i_1}, \ldots, e_{i_k}, \ldots, f_m(e_{i_m}), \ldots, e_{i_n}],$$
то есть, 
$$[f_k(e_i) -f_m(e_i), e_{i_1}, \ldots, e_{i_k}, \ldots, \hat{e}_{i_m}, \ldots, e_{i_n}]=0.$$
Изменяя  индексы $k$ и $m$, а также множество $I=\{i_l\} \subset \{1, \ldots, n+1\},$ где $|I|=n-1,$
мы можем получить, что для произвольных $k$ и $m$, а также для $j \neq i$ верно $$(f_k(e_i) -f_m(e_i))|_{e_j}=0.$$

Пусть $A=(a_{ij})$ --- матрица преобразования $f_0$ и $B_k=(b_{ij}^k)$ --- матрица преобразования $f_k$. По доказанному выше, 
мы можем считать, что $b_{ij}^k=b_{ij}$ при $i \neq j$. Тогда,
$$(-1)^{n+i+1}f_0(e_i)=f_0[e_1, \ldots, \hat{e}_i, \ldots, e_{n+1}]=$$
$$\sum\limits_{j}[e_1, \ldots,f_j( e_j), \ldots, \hat{e}_i, \ldots, e_{n+1}]=
(-1)^{n+i+1}\left(\sum\limits_{k,k \neq i} b_{kk}^k\right) e_i  + \sum\limits_{j,j \neq i} (-1)^{n+i}b_{ji} e_j,$$
то есть 
$a_{ij}=-b_{ji},$ при $i \neq j.$
Пусть $A^*$ матрица, составленная из элементов матрицы $A$ с нулевыми элементами по диагонали и $g_0$ линейное отображение с матрицей $A^*.$
Соответственно, $B^*$ матрица, составленная из элементов матрицы $B_k$, 
с нулевыми элементами по диагонали и $g$ линейное отображение с матрицей $B^*.$
Согласно лемме 1, $(g_0, g, \ldots, g)$ --- $(n+1)$-арное дифференцирование алгебры $A_{n+1}$. Таким образом, 
разность $$(f_0, f_1, \ldots, f_n) - (g_0, g, \ldots, g)=(f_0^*, f_1^*, \ldots, f_n^*)$$
будет являться $(n+1)$-арным дифференцированием, где каждый элемент $f_k^*$ в базисе $\{e\}$ имеет диагональную матрицу линейного преобразования, 
то есть, удовлетворяет условию леммы.
Лемма доказана.

\medskip

{\bf  Лемма 3.}  В терминах леммы 2, верно 
$$(f_0^*, f_1^*, \ldots, f_n^*)=(\sum h_j \cdot id, h_1 \cdot id, \ldots, h_n \cdot id)+ (d_0, d, \ldots, d).$$

\medskip

{\bf Доказательство.} 
Будем считать, что $f_i^*(e_j)=f_i^j e_j.$ 
Ясно, что $$(\sum f_i^{n+1} \cdot id, f_1^{n+1} \cdot id, \ldots, f_n^{n+1} \cdot id)$$ является $(n+1)$-арным дифференцированием, которое мы будем обозначать
$(\sum h_j \cdot id, h_1 \cdot id, \ldots, h_n \cdot id)$.
Рассмотрим разность $(n+1)$-арных дифференцирований $(f_0^*, f_1^*, \ldots, f_n^*)$ и $(\sum h_j \cdot id, h_1 \cdot id, \ldots, h_n \cdot id).$
Полученное $(n+1)$-арное дифференцирование обозначим $(d_0, d_1, \ldots, d_n)$, где $d_i(e_j)=d_i^j e_j$ и $d_i(e_{n+1})=0.$ 
Для наглядности мы можем записать коэффициенты $d_i^j$  в виде матрицы 
\begin{eqnarray}
\label{ma} \left(\begin{array}{cccrcccc}
d_1^1 	   & d_1^2 &  \ldots 			& d_{1}^{n-2}	& d_{1}^{n-1} & d_1^n 	&|& d_1^{n+1}=0  \\
d_2^1 	   & d_2^2 &   			& d_{2}^{n-2}	& d_{2}^{n-1} & d_2^n 		&|& d_2^{n+1}=0 \\

\vdots	  &  	& \ddots    	&  	& \vdots               & &               &\\
d_{n-3}^1 & d_{n-3}^2 	&  		&d_{n-3}^{n-2}	& d_{n-3}^{n-1} & d_{n-3}^n 	&|& d_{n-3}^{n+1}=0  \\
d_{n-2}^1 & d_{n-2}^2 	&  \ldots 	&d_{n-2}^{n-2}	& d_{n-2}^{n-1} & d_{n-2}^n 	&|& d_{n-2}^{n+1}=0  \\
d_{n-1}^1 & d_{n-1}^2 	&  		& d_{n-1}^{n-2}	& d_{n-1}^{n-1} & d_{n-1}^n 	&|& d_{n-1}^{n+1}=0  \\
d_n^1 	  & d_n^2 	&    		&d_{n}^{n-2}	& d_{n}^{n-1} 	& d_n^n 	&|& d_n^{n+1}=0 
\end{array} \right).
\end{eqnarray}
Отметим, что 
$$0=d_0[e_1, \ldots, e_k, \ldots, e_m, \ldots, e_n]+d_0[e_1, \ldots, e_m, \ldots, e_k, \ldots, e_n]=$$
$$(d_k^k+d_m^m-d_m^k-d_k^m)e_{n+1},$$
то есть 
\begin{eqnarray}
\label{ss1} d_k^k+d_m^m=d_k^m+d_m^k, \mbox{ где }k,m \leq  n.\end{eqnarray}
Аналогично мы можем получить
$$0=d_0[e_2, \ldots, e_k, \ldots, e_m, \ldots, e_{n+1}]+d_0[e_2, \ldots, e_m, \ldots, e_k, \ldots, e_{n+1}]=$$
$$(-1)^n(d_{m-1}^k+d_{k-1}^m-d_{m-1}^m-d_{k-1}^k)e_{1},$$
то есть,
\begin{eqnarray}
\label{ss2} d_{m-1}^k+d_{k-1}^m=d_{m-1}^m+d_{k-1}^k, \mbox{ где } m,k \geq 2.
\end{eqnarray}

Ясно, что значения в выражениях (\ref{ss1}-\ref{ss2}) образуют вершины квадрата
в матрице (\ref{ma}) с диагональю лежащей на диагонали матрицы, либо на побочной диагонали.

Осталось показать, что $d_1=\ldots=d_n.$ Для этого необходимо показать, что 
для произвольного $j$ верно $d_i^j=d_k^j$. 
Отметим, что в силу (\ref{ss2}) мы имеем 
$$d_{m-1}^m=d_n^m, m \geq 2.$$

Откуда, по (\ref{ss1}) и $d_{n-1}^n=d_n^n$, мы получим $d_n^{n-1}=d_{n-1}^{n-1}.$
Теперь заметим, что $$d_{n-2}^{n-1}=d_{n-2}^{n-2}$$ и верно условие (\ref{ss2}), следовательно $$d_{n-1}^{n-1}=d_{n-1}^{n-2}.$$
Отсюда и (\ref{ss2}), видим $$d^{n}_{n-2}=d^{n}_{n-1}=d^{n}_{n}.$$
Последнее, с учетом (\ref{ss1}), влечет $$d_{n-2}^{n-2}=d_{n}^{n-2}.$$ 
Продолжая дейстововать аналогичным образом, мы получим, что для произвольного $j$ верно $d_i^j=d_k^j$, 
то есть требуемое. Лемма доказана.

\medskip

\section{Основные результаты.}

Теперь заметим, что если $(f_0,f, \ldots, f)$ --- $(n+1)$-арное дифференцирование $n$-арной алгебры $A_{n+1},$
где каждое из отображений $f_0$ и $f$ на базисных элементах действует как $f_0(e_i)=f_0^i e_i, f(e_i)=f_i e_i$,
то $$f_0^i e_i= (-1)^{n+1+i}f_0[e_1, \ldots, \hat{e}_i, \ldots, e_{n+1}]=(\sum\limits_{j=1}^{n+1} f_j - f_i)e_i.$$
Откуда получаем, что 
$$f_0^i = \sum\limits_j f_j - f_i,$$
и, следовательно, 
$$\sum\limits_j f_0^j = n \sum\limits_j f_j.$$
Полученные результаты, дают 
\begin{eqnarray}\label{osnre}
f_i=\frac{1}{n}\sum\limits_j f_0^j - f_0^i. 
\end{eqnarray}
То есть, отображение $f$ однозначно определяется из вида отображения $f_0$ и произвольное отображение
$f_0$, такое что $f_0(e_i)=f_0^ie_i$ задает $(n+1)$-арное дифференцирование $(f_0, f, \ldots, f)$, где 
$f_0$ и $f$ связаны посредством соотношения (\ref{osnre}). 
Таким образом, исходя из полученного и леммы 1, мы имеем 

\medskip

{\bf  Теорема 4.}  Пусть $A$ --- простая конечномерная алгебра Филиппова над алгебраически замкнутым полем характеристики нуль, 
тогда $\Delta(A) = sl_{n+1}.$

\medskip

Полученные результаты мы можем подытожить в следующей теореме.

\medskip

{\bf  Теорема 5.}  Пусть $(f_0, f_1, \ldots, f_n)$ --- $(n+1)$-арное дифференцирование простой 
конечномерной $n$-арной алгебры Филиппова над алгебраически замкнутым полем характеристики нуль, 
тогда 
$$(f_0, f_1, \ldots, f_n)= 
(\sum\limits_{j=1}^{n} h_j \cdot id , h_1 \cdot id, \ldots, h_n \cdot id)+ (d_0, d, \ldots, d)$$
и $
[d_0]^{T}+[d]=0$, где $[d_0],[d]$ --- матрицы линейных отображений $d_0,d.$



\medskip 

Отметим, что из теоремы 5 легко получается описание $\delta$-дифференцирований простых конечномерных алгебр Филиппова 
над алгебраически замкнутым полем характеристики нуль, 
которые согласуются с результатами \cite{kay_nary}.

\medskip

{\bf  Теорема 6.}  
Пусть $(f_0, f_1, \ldots, f_n)$ --- $(n+1)$-арное дифференцирование полупростой 
конечномерной $n$-арной алгебры Филиппова $A$ над алгебраически замкнутым полем характеристики нуль, 
тогда $A=\oplus_{i=1}^{t} I_i,$ где $I_i$ --- простой идеал алгебры $A$
и 

1) $f_j(I_i) \subseteq I_i;$

2) матрицы линейных отображений $f_j$ в подходящем базисе имеют блочно-диагональный вид, 
где каждый квадратный блок имеет размерность $n+1$ и строится исходя из описания 
соответствующего $(n+1)$-арного дифференцирования простой алгебры Филиппова $I_i$ посредством теоремы 5,
то есть в подходящем базисе
$$(f_0, f_1, \ldots, f_n)= 
(\sum\limits_{j=1}^{n} h_j, h_1, \ldots, h_n)+ (d_0, d, \ldots, d),$$
где 
$$[h_i]=
 \left(\begin{array}{ccrc}
h_i^1 E_{n+1} 	    & \ldots 			& 0	 \\
\vdots	  	    	& \ddots    	  	& \vdots                              \\
0 	  	    &   \ldots 			 & h_i^t E_{n+1}   \\
\end{array} \right) \mbox{ и }
[d_0]^{T}+[d]=0.$$
Здесь под $E_{n+1}$ подразумевается единичная матрица размера $n+1$, а через $[P]$ 
обозначена матрица линейного преобразования $P$.


\medskip

{\bf Доказательство.} 
Согласно результатам \cite{Ling}, $A=\oplus_{i=1}^{t} I_i,$ где $I_i$ --- простой идеал алгебры $A$ и $I_i\cong A_{n+1}.$

Пусть $A=I \oplus J,$ где $J$ --- простой идеал, а $I$ --- прямая сумма некоторого количества идеалов изоморфных $J$. 
Тогда для элементов $i_k \in I$ верно 
$$f_0[i_1, \ldots , i_n] =\sum[i_1, \ldots, f_k(i_k), \ldots, i_n] \in I.$$
Учитывая классификацию простых конечномерных алгебр Филиппова \cite{Ling} и структуру алгебры $A_{n+1}$,
мы можем заключить, что $f_0(I) \subseteq I.$
Допустим, что для некоторого $t$ верно $f_t(I)|_J \neq 0.$
Тогда 
$$0=f_0[J, \ldots, J, I, J, \ldots, J]=[J, \ldots, J, f_t(I), J, \ldots, J],$$
то есть, $f_t(I)|_J$ идеал в $J.$ 
Откуда следует, что либо $f_t(I) \subseteq I,$ либо $f_t(I)|_J=J.$
В силу того, что алгебра $J$ является простой и, следовательно, 
не является нильпотентной, то второй случай не возможен. 
Таким образом, мы показали, что выполнено условие 1. 

Для доказательства второго условия теоремы 
выберем базис алгебры $A$ как объединение базисов $I_i$.
Нам достаточно рассмотреть ограничение отображений $f_t$ на $I_i$ и воспользоваться теоремой 5.
Теорема доказана.

\medskip

Из теорем 4 и 6 легко получается

\medskip

{\bf  Теорема 7.}  Пусть $A$ --- полупростая конечномерная $n$-арная алгебра Филиппова над алгебраически замкнутым полем характеристики нуль, 
тогда $\Delta(A) = \oplus sl_{n+1}.$

\medskip

Учитывая результаты теорем 5 и 7, а также \cite{LL}, мы можем сформулировать следующую гипотезу.

\medskip

{\bf  Гипотеза.} Если $n$-арная алгебра Филиппова обладает свойством $GDer(A)=End(A),$
то либо ее размерность не выше чем $n$, либо она простая $(n+1)$-мерная алгебра $A_{n+1}$.
 
\medskip

В заключение, 
автор выражает благодарность проф. В. Н. Желябину, проф. А. П. Пожидаеву и проф. П. С. Колесникову 
за внимание к работе и конструктивные замечания.


\end{document}